\documentclass{article}

\usepackage[letterpaper]{geometry}
\usepackage{amsmath,bussproofs,float}

\floatstyle{boxed}
\restylefloat{figure}

\newtheorem{theorem}{Theorem}[section]

\renewcommand{\fCenter}{\Rightarrow}
\EnableBpAbbreviations

\newcommand{\AXM}[1]{\AXC{$#1$}}

\newcommand{\AXN}[1]{\AX$#1$}
\newcommand{\BIN}[1]{\BI$#1$}
\newcommand{\UIN}[1]{\UI$#1$}
\newcommand{\TIN}[1]{\TI$#1$}

\newcommand{\ph}{\varphi}
\newcommand{\ex}[1]{\exists #1 \;} 
\newcommand{\fa}[1]{\forall #1 \;} 
\newcommand{\na}[1]{\mathit{#1}}    

\title{Proof Theory}
\author{Jeremy Avigad}

\begin{document}

\maketitle

\begin{abstract}
Proof theory began in the 1920's as a part of Hilbert's program, which aimed to secure the foundations of mathematics by modeling infinitary mathematics with formal axiomatic systems and proving those systems consistent using restricted, finitary means. The program thus viewed mathematics as a system of reasoning with precise linguistic norms, governed by rules that can be described and studied in concrete terms. Today such a viewpoint has applications in mathematics, computer science, and the philosophy of mathematics. 

\medskip \noindent \emph{Keywords:} proof theory, Hilbert's program, foundations of mathematics
\end{abstract}

\section{Introduction}
\label{introduction:section}

At the turn of the nineteenth century, mathematics exhibited a style of argumentation that was more explicitly computational than is common today. Over the course of the century, the introduction of abstract algebraic methods helped unify developments in analysis, number theory, geometry, and the theory of equations, and work by mathematicians like Richard Dedekind, Georg Cantor, and David Hilbert towards the end of the century introduced set-theoretic language and infinitary methods that served to downplay or suppress computational content. This shift in emphasis away from calculation gave rise to concerns as to whether such methods were meaningful and appropriate in mathematics. The discovery of paradoxes stemming from overly naive use of set-theoretic language and methods led to even more pressing concerns as to whether the modern methods were even consistent. This led to heated debates in the early twentieth century and what is sometimes called the ``crisis of foundations.''

In lectures presented in 1922, Hilbert launched his \emph{Beweistheorie}, or Proof Theory, which aimed to justify the use of modern methods and settle the problem of foundations once and for all. This, Hilbert argued, could be achieved as follows:
\begin{itemize}
\item First, represent portions of the abstract, infinitary mathematical reasoning in question using formal axiomatic systems, which prescribe a fixed formal language and precise rules of inference.
\item Then view proofs in these systems as finite, combinatorial objects, and prove the consistency of such systems---i.e.~the fact that there is no way to derive a contradiction---using unobjectionable, concrete arguments.
\end{itemize}
In doing so, said Hilbert,
\begin{quote}
 \ldots we move to a higher level of contemplation, from which the axioms, formulae, and proofs of the mathematical theory are themselves the objects of a contentional investigation. But for this purpose the usual contentual ideas of the mathematical theory must be replaced by formulae and rules, and imitated by formalisms. In other words, we need to have a strict formalization of the entire mathematical theory\ldots. In this way the contentual thoughts (which of course we can never wholly do without or eliminate) are removed elsewhere---to a higher plane, as it were; and the same time it becomes possible to draw a sharp and systematic distinction in mathematics between the formulae and formal proofs on the one hand, and the contentual ideas on the other. \cite{hilbert:22}
\end{quote}
G\"odel's second incompleteness theorem shows that any ``unobjectionable'' portion of mathematics is insufficient to establish its own consistency, let alone the consistency of any theory properly extending it. Although this dealt a blow to Hilbert's program as it was originally formulated, the more general project of studying mathematical reasoning in syntactic terms, especially with respect to questions of algorithmic or otherwise concrete content, has been fruitful. Moreover, the general strategy of separating syntactic and semantic concerns and of maintaining a syntactic viewpoint where possible has become a powerful tool in formal epistemology. (See \cite{sieg:99,zach:06} for more on Hilbert's program.)

Today, Proof Theory can be viewed as the general study of formal deductive systems. Given that formal systems can be used to model a wide range of types of inference---modal, temporal, probabilistic, inductive, defeasible, deontic, and so on---work in the field is varied and diverse. Here I will focus specifically on the proof theory of \emph{mathematical} reasoning, but even with this restriction, the field is dauntingly broad: the 1998 \emph{Handbook of Proof Theory} \cite{buss:98} runs more than 800 pages, with a name index that is almost as long as this article. As a result, I can only attempt to convey a feel for the subject's goals and methods of analysis, and help ease the reader into the broader literature. References are generally to surveys and textbooks, and results are given without attribution. 

In Section~\ref{deduction:section}, I describe natural deduction and a sequent calculus for first-order logic, and state the cut-elimination theorem and some of its consequences. This is one of the field's most fundamental results, and provides a concrete example of proof-theoretic method. In Section~\ref{tools:section}, I survey various aspects of proof-theoretic analysis, and, finally, in Section~\ref{applications:section}, I discuss some applications.

\section{Natural deduction and sequent calculi}
\label{deduction:section}

I will assume the reader is familiar with the language of first-order logic. Contemporary logic textbooks often present formal calculi for first-order logic with a long list of axioms and a few simple rules, but these are generally not very convenient for modeling deductive arguments or studying their properties. A system which fares better on both counts in given by Gerhard Gentzen's system of \emph{natural deduction}, a variant of which we will now consider.

Natural deduction is based on two fundamental observations. The first is that it is natural to describe the meaning, or appropriate use, of a logical connective by giving the conditions under which one can \emph{introduce} it, that is, derive a statement in which that connective occurs, and the methods by which one can \emph{eliminate} it, that is, draw conclusions from statements in which it occurs. For example, one can establish a conjunction $\ph \wedge \psi$ by establishing both $\ph$ and $\psi$, and, conversely, if one assumes or has previously established $\ph \wedge \psi$, one can conclude either $\ph$ or $\psi$, at will. 

The second observation is that it is natural to model logical arguments as taking place under the context of a list of hypotheses, either implicit or explicitly stated. If $\Gamma$ is a finite set of hypotheses and $\ph$ is a first-order formula, the \emph{sequent} $\Gamma \fCenter \ph$ is intended to denote that $\ph$ follows from $\Gamma$. For the most part, these hypotheses stay fixed over the course of an argument, but under certain circumstances they can be removed, or \emph{canceled}. For example, one typically proves an implication $\ph \rightarrow \psi$ by temporarily assuming that $\ph$ holds and arguing that $\psi$ follows. The introduction rule for implication thus reflects the fact that deriving $\psi$ from a set of hypotheses $\Gamma$ together with $\ph$ is the same as deriving $\ph \rightarrow \psi$ from $\Gamma$. 

Writing $\Gamma, \ph$ as an abbreviation for $\Gamma \cup \{ \ph \}$, the rules for natural deduction are shown in Figure 1.
\begin{figure}[t]
\begin{center}
\begin{tabular}{cc}
\\
\multicolumn{2}{c}{
\AXM{}
\UIN{\Gamma, \ph \fCenter \ph}
\DP
}
\\
\\
\AXN{\Gamma \fCenter \ph}
\AXN{\Gamma \fCenter \psi}
\BIN{\Gamma \fCenter \ph \wedge \psi}
\DP
&
\AXN{\Gamma \fCenter \ph_0 \wedge \ph_1}
\UIN{\Gamma \fCenter \ph_i}
\DP
\\
\\
\AXN{\Gamma \fCenter \ph_i}
\UIN{\Gamma \fCenter \ph_0 \lor \ph_1}
\DP
&
\AXN{\Gamma \fCenter \ph \lor \psi}
\AXN{\Gamma, \ph \fCenter \theta}
\AXN{\Gamma, \psi \fCenter \theta}
\TIN{\Gamma \fCenter \theta}
\DP
\\
\\
\AXN{\Gamma, \ph \fCenter \psi}
\UIN{\Gamma \fCenter \ph \rightarrow \psi}
\DP
&
\AXN{\Gamma \fCenter \ph \rightarrow \psi}
\AXN{\Gamma \fCenter \ph}
\BIN{\Gamma \fCenter \psi}
\DP
\\
\\
\AXN{\Gamma \fCenter \ph}
\UIN{\Gamma \fCenter \fa y \ph[y/x]}
\DP
&
\AXN{\Gamma \fCenter \fa x \ph}
\UIN{\Gamma \fCenter \ph[t/x]}
\DP
\\
\\
\AXN{\Gamma \fCenter \ph[t/x]}
\UIN{\Gamma \fCenter \ex x \ph}
\DP
&
\AXN{\Gamma \fCenter \ex y \ph[y/x]}
\AXN{\Gamma, \ph \fCenter \psi}
\BIN{\Gamma \fCenter \psi}
\DP
\end{tabular}
\end{center}
\caption{Natural deduction. Derivability of a sequent $\Gamma \Rightarrow \ph$ means that $\ph$ is a consequence of the set of hypotheses $\Gamma$, and $\Gamma, \ph$ denotes $\Gamma \cup \{ \ph \}$.}
\end{figure}
The quantifier rules are subject to the usual restrictions. For example, in the introduction rule for the universal quantifier, the variable $x$ cannot be free in any hypothesis. For intuitionistic logic, one also needs the rule {\em ex falso sequitur quodlibet}, which allows one to conclude $\Gamma \fCenter \ph$ from $\Gamma \fCenter \bot$, where $\bot$ represents falsity. One can then define negation, $\lnot \ph$, as $\ph \rightarrow \bot$. For classical logic, one adds {\em reductio ad absurdum}, or proof by contradiction, which allows one to conclude $\Gamma \fCenter \ph$ from $\Gamma, \lnot \ph \fCenter \bot$.


For many purposes, however, \emph{sequent calculi} provide a more convenient representation of logical derivations. Here, sequents are of the form $\Gamma \fCenter \Delta$, where $\Gamma$ and $\Delta$ are finite sets of formulas, with the intended meaning that the conjunction of the hypotheses in $\Gamma$ implies the \emph{disjunction} of the assertions in $\Delta$. The rules are as shown in Figure 2.
\begin{figure}[t]
\begin{center}
\begin{tabular}{cc}
\multicolumn{2}{c}{
\AXM{}
\UIN{\Gamma, \ph \fCenter \Delta, \ph} 
\DP} \\
\\
\AXN{\Gamma, \ph_i \fCenter \Delta} 
\UIN{\Gamma, \ph_0 \wedge \ph_1 \fCenter \Delta}
\DP &
\AXN{\Gamma \fCenter \Delta, \ph}
\AXN{\Gamma \fCenter \Delta, \psi}
\BIN{\Gamma \fCenter \Delta, \ph \wedge \psi}
\DP \\
\\
\AXN{\Gamma, \ph \fCenter \Delta}
\AXN{\Gamma, \theta \fCenter \Delta}
\BIN{\Gamma, \ph \lor \theta \fCenter \Delta}
\DP &
\AXN{\Gamma \fCenter \Delta, \ph_i}
\UIN{\Gamma \fCenter \Delta, \ph_0 \lor \ph_1}
\DP \\
\\
\AXN{\Gamma, \fCenter \Delta, \ph}
\AXN{\Gamma, \theta \fCenter \Delta}
\BIN{\Gamma, \ph \rightarrow \theta \fCenter \Delta}
\DP &
\AXN{\Gamma, \ph \fCenter \Delta, \psi}
\UIN{\Gamma \fCenter \Delta, \ph \rightarrow \psi}
\DP \\
\\
\AXN{\Gamma, \ph[t/x] \fCenter \Delta}
\UIN{\Gamma, \fa x \ph  \fCenter \Delta}
\DP &
\AXN{\Gamma \fCenter \Delta, \psi[y/x]}
\UIN{\Gamma \fCenter \Delta, \fa x \psi}
\DP \\
\\
\AXN{\Gamma, \ph[y/x] \fCenter \Delta}
\UIN{\Gamma, \ex x \ph \fCenter \Delta}
\DP &
\AXN{\Gamma \fCenter \Delta, \psi[t/x]}
\UIN{\Gamma \fCenter \Delta, \ex x \psi}
\DP \\
\\
\multicolumn{2}{c}{
\AXN{\Gamma \fCenter \Delta, \ph}
\AXN{\Gamma, \ph \fCenter \Delta}
\BIN{\Gamma \fCenter \Delta}
\DP}
\end{tabular}
\end{center}
\caption{The sequent calculus.}
\end{figure}
The last rule is called the \emph{cut rule}: it is the only rule containing a formula in the hypothesis that may be entirely unrelated to the formulas in the conclusion. Proofs that do not use the cut rule are said to be \emph{cut free}. One obtains a proof system for intuitionistic logic by restricting $\Delta$ to contain at most one formula, and adding an axiomatic version of \emph{ex falso sequitur quodlibet}: $\Gamma, \bot \fCenter \ph$. The cut-elimination theorem is as follows:

\begin{theorem}
 If $\Gamma \fCenter \Delta$ is derivable in the sequent calculus with cut, then it is derivable without cut. 
\end{theorem}

Gentzen's proof gives an explicit algorithm for removing cuts from a proof. The algorithm, unfortunately, can yield an iterated exponential increase in the size of proofs, and one can show that there are cases in which such an increase cannot be avoided. The advantage of having a cut-free proof is that the formulas in each sequent are built up directly from the formulas in the sequents above it, making it easy to extract useful information. For example, the following are two consequences of the cut-elimination theorem, easily proved by induction on cut-free proofs. 

The first is known as \emph{Herbrand's theorem}. Recall that a formula of first-order logic is said to be \emph{existential} if it consists of a block of existential quantifiers followed by a quantifier-free formula. Similarly, a formula is said to be \emph{universal} if it consists of a block of universal quantifiers followed by a quantifier-free formula. Herbrand's theorem says that if it is possible to prove an existential statement from some universal hypotheses, then in fact there is an explicit sequence of terms in the language that witness the truth of the conclusion.

\begin{theorem}
\label{herbrand:thm}
 Suppose $\ex {\vec x} \ph(\vec x)$ is derivable in classical first-order logic from a set of hypotheses $\Gamma$, where $\ph$ is quantifier-free and the sentences in $\Gamma$ are universal sentences. Then there are sequences of terms $\vec t_1, \vec t_2, \ldots, \vec t_k$ such that the disjunction $\ph(\vec t_1) \lor \ph(\vec t_2) \lor \ldots \lor \ph(\vec t_k)$ has a quantifier-free proof from instances of the sentences in $\Gamma$.
\end{theorem}

For intuitionistic logic, one has a stronger property, known as the \emph{explicit definability property}.

\begin{theorem}
\label{ed:thm}
  Suppose $\ex {\vec x} \ph(\vec x)$ is derivable in intuitionistic first-order logic from a set of hypotheses $\Gamma$ in which neither $\lor$ nor $\exists$ occurs in a strictly positive part. Then there are terms $\vec t$ such that $\ph(\vec t)$ is also derivable from $\Gamma$.
\end{theorem}

Theorem~\ref{herbrand:thm} provides a sense in which explicit information can be extracted from certain classical proofs, and Theorem~\ref{ed:thm} provides a sense in which intuitionistic logic is constructive. We have thus already encountered some of the central themes of proof-theoretic analysis: 
\begin{itemize}
\item Important fragments of mathematical reasoning can be captured by formal systems.
\item One can study the properties of these formal systems, for example, describing transformations of formulas and proofs, translations between formulas and proofs in different systems, and canonical normal forms for formulas and proofs.
\item The methods provide information about the logic that is independent of the choice of formal system that is used to represent it.
\end{itemize}

For more on the cut-elimination theorems, see \cite{buss:98b,schwichtenberg:77,takeuti:87,troelstra:schwichtenberg:00}. 

\section{Methods and goals}
\label{tools:section}

\subsection{Classical foundations}

Recall that Hilbert's program, broadly construed, involves representing mathematical reasoning in formal systems and then studying those formal systems as mathematical objects themselves. The first step, then, requires finding the right formal systems. It is common today to view mathematical reasoning as consisting of a properly mathematical part that is used in conjunction with more general forms of logical reasoning, though there are still debates as to where to draw the line between the two. In any case, the following list portrays some natural systems of reasoning in increasing logical/mathematical strength:
\begin{enumerate}
\item pure first-order logic
\item primitive recursive arithmetic (denoted $\na{PRA}$)
\item first-order arithmetic ($\na{PA}$)
\item second-order arithmetic ($\na{PA^2}$)
\item higher-order arithmetic ($\na{PA^\omega}$)
\item Zermelo-Fraenkel set theory ($\na{ZF}$)
\end{enumerate}

Primitive recursive arithmetic was designed by Hilbert and Bernays to be a patently finitary system of reasoning. The system allows one to define functions on the natural numbers using a simple schema of primitive recursion, and prove facts about them using a principle of induction:
\[
\ph(0) \wedge \fa x (\ph(x) \rightarrow \ph(x + 1)) \rightarrow \fa x \ph(x).
\]
In words, if $\ph$ holds of $0$ and, whenever it holds of some number, $x$, it holds of $x + 1$, then $\ph$ holds of every number. Here $\ph$ is assumed to be a quantifier-free formula. In fact, one can replace this axiom with a suitable induction \emph{rule}, whereby primitive recursive arithmetic can be formulated without quantifiers at all. Surprisingly, via coding of finitary objects as natural numbers, this system is expressive and strong enough to develop most portions of mathematics that involve only finite objects and structures \cite{avigad:03}. Peano arithmetic can be viewed as the extension of $\na{PRA}$ with induction for all first-order formulas.

There is no effective axiomatization of second- or higher-order logic that is complete for the standard semantics (where, for example, second-order quantifiers are assumed to range over all subsets of the universe of individuals). As a result, one has to distinguish axiomatic second- and higher-order logic from the corresponding semantic characterization. Axiomatically, one typically augments first-order logic with comprehension rules that assert that every formula defines a set (or predicate):
\[
\ex X \fa y (X (y) \leftrightarrow \ph)
\]
Here $\ph$ is a formula in which $X$ does not occur, although $\ph$ is allowed to have other free variables in addition to $y$. One can augment these with suitable choice principles as well. Second-order arithmetic can be viewed as the extension of Peano arithmetic with second-order logic and second-order principles of induction, but one can, alternatively, interpret second-order arithmetic in second-order logic together with an axiom asserting the existence of an infinite domain. Similar considerations hold for higher-order logic as well. 

Axioms for set theory can be found in any introductory set theory textbook, such as \cite{kunen:80}. Of course, these axioms can be extended with stronger hypotheses, such as large cardinal axioms. For information on primitive recursive arithmetic, see \cite{goodstein:57,troelstra:schwichtenberg:00}; for first-order arithmetic, see \cite{buss:98b,hajek:pudlak:93,kaye:91}; for second-order arithmetic, see \cite{simpson:99}; for higher-order arithmetic, see \cite{takeuti:87}.

\subsection{Constructive foundations}
\label{type:theory}

Given the history of Hilbert's program, it should not be surprising that proof theorists have also had a strong interest in formal representations of constructive and intuitionistic reasoning. From an intuitionistic standpoint, the use of the excluded middle, $\ph \vee \lnot \ph$, is not acceptable, since, generally speaking, one may not know (or have an algorithm to determine) which disjunct holds. For example, in classical first-order arithmetic, one is allowed to assert $\ph \vee \lnot \ph$ for a formula $\ph$ that expresses the twin primes conjecture, even though we do not know which is the case. If one restricts the underlying logic to intuitionistic logic, however, one obtains \emph{Heyting arithmetic}, which is constructively valid. 

Stronger systems tend to be based on what has come to be known as the Curry-Howard-Tait \emph{propositions as types} correspondence. The idea is that, from a constructive perspective, any proposition can be viewed as specifying a type of data, namely, the type of construction that warrants the claim that the proposition is true. A proof of the proposition is thus a construction of the corresponding type. For example, a proof of $\ph \wedge \psi$ is a proof of $\ph$ paired with a proof of $\psi$, and so $\ph \wedge \psi$ corresponds to the type of data consisting of pairs of type $\ph$ and $\psi$. Similarly, a proof of $\ph \rightarrow \psi$ should be a procedure transforming a proof of $\ph$ into a proof of $\psi$, so $\ph \rightarrow \psi$ corresponds to a type of functions. This gives rises to systems of \emph{constructive type theory,} of which the most important examples are \emph{Martin-L\"of type theory} and an impredicative variant designed by Coquand and Huet, the \emph{calculus of constructions}. Thus, our representative sample of constructive proof systems, in increasing strength, runs as follows:
\begin{enumerate}
\item intuitionistic first-order logic
\item primitive recursive arithmetic ($\na{PRA}$)
\item Heyting arithmetic ($\na{HA}$)
\item Martin-L\"of type theory ($\na{ML}$)
\item the calculus of inductive constructions ($\na{CIC}$)
\end{enumerate}
Good references for intuitionistic systems in general are \cite{beeson:85,troelstra:van:dalen:88both}. For more information on type theory, see \cite{sambin:98}; for the calculus of inductive constructions in particular, see \cite{bertot:casteran:04}.

\subsection{Reverse mathematics}

In the 1970's, Harvey Friedman observed that by restricting the induction and comprehension principles in full axiomatic second-order arithmetic, one obtains theories that are strong enough, on the one hand, to represent significant parts of ordinary mathematics, but weak enough, on the other hand, to be amenable to proof-theoretic analysis. He then suggested calibrating various mathematical theorems in terms of their axiomatic strength. Whereas in ordinary (meta)mathematics, one proves theorems from axioms, Friedman noticed that it is often the case that a mathematical theorem can be used in the other direction, namely, to prove an underlying set-existence principle, over a weak base theory. That is, it is often the case that a theorem of mathematics is formally \emph{equivalent} to a set comprehension principle that is used to prove it.

In that years that followed, Friedman, Stephen Simpson, and many others worked to calibrate the axiomatic assumptions used in a wide range of subjects. They isolated five key theories along the way: 
\begin{enumerate}
\item $\na{RCA_0}$: a weak base theory, conservative over primitive recursive arithmetic, with a \emph{recursive comprehension axiom}, that is, a principle of comprehension for recursive (computable) sets. 
\item $\na{WKL_0}$: adds \emph{weak K\"onig's lemma,} a compactness principle, to $\na{RCA_0}$. 
\item $\na{ACA_0}$: adds the \emph{arithmetic comprehension axiom}, that is, comprehension for arithmetically definable sets. 
\item $\na{ATR_0}$: adds a principle of \emph{arithmetical transfinite recursion,} which allows one to iterate arithmetic comprehension along countable well-orderings. 
\item $\na{\Pi^1_1\mathord{-}CA}$: adds the $\Pi^1_1$ \emph{comprehension axiom}, that is, comprehension for $\Pi^1_1$ sets. 
\end{enumerate} 
Simpson \cite{simpson:99} provides the best introduction to these theories and the reverse mathematics program.

\subsection{Comparative analysis and reduction}

We have now seen a sampling of the many formal systems that have been designed to formalize various aspects of mathematics. Proof theorists have also invested a good deal of energy in understanding the relationships between the systems. Often, results take the form of \emph{conservation theorems} which fit the following pattern, where $T_1$ and $T_2$ are theories and $\Gamma$ is a class of sentences:
\begin{quote}
Suppose $T_1$ proves a sentence $\ph$, where $\ph$ is in $\Gamma$. Then $T_2$ proves it as well (or perhaps a certain translation, $\ph'$).
\end{quote}
Such a result, when proved in a suitably restricted base theory, provides a foundational reduction of the theory $T_1$ to $T_2$, justifying the principles of $T_1$ relative to $T_2$. For example, such theorems can be used to reduce:
\begin{itemize}
\item an infinitary theory to a finitary one
\item a nonconstructive theory to a constructive one
\item an impredicative theory to a predicative one
\item a nonstandard theory (in the sense of nonstandard analysis) to a standard one
\end{itemize}
For example:
\begin{enumerate}
\item Versions of primitive recursive arithmetic based on classical, intuitionistic, or quantifier-free logic, all prove the same $\Pi_2$ theorems (in an appropriate sense) \cite{troelstra:schwichtenberg:00}.
\item The G\"odel-Gentzen double-negation interpretation and variations, like the Friedman-Dragalin A-translation, interpret a number of classical systems in intuitionistic ones, like $\na{PA}$ in $\na{HA}$ \cite{avigad:00b, buss:98e,feferman:77,troelstra:schwichtenberg:00,troelstra:van:dalen:88both}.
\item There are various translations between theories in the language of (first-, second-, or higher-order) arithmetic and subsystems of set theory \cite{pohlers:09,simpson:99}.
\item Both $\na{I\Sigma_1}$, the subsystem of Peano arithmetic in which induction is restricted to $\Sigma_1$ formulas, and $\na{WKL_0}$, the subsystem of second-order arithmetic based on Weak K\"onig's Lemma, are conservative over primitve recursive arithmetic for the class of $\Pi_2$ sentences \cite{avigad:02c,buss:98b,hajek:pudlak:93,kohlenbach:08,sieg:85,simpson:99,troelstra:schwichtenberg:00}.
\item Cut elimination or an easy model-theoretic argument shows that a restricted second-order version, $\na{ACA_0}$, of Peano arithmetic is a conservative extension of Peano arithmetic itself. Similarly, G\"odel-Bernays-von Neumann set theory $\na{GBN}$, which has both sets and classes, is a conservative extension of Zermelo-Fraenkel set theory. See, for example, \cite{pudlak:98,simpson:99}. In general, proofs in $\na{ACA_0}$ may suffer an iterated exponential increase in length when translated to $\na{PA}$, and similarly for $\na{GBN}$ and $\na{ZF}$, or $\na{I\Sigma_1}$ and $\na{PRA}$.
\item Theories of nonstandard arithmetic and analysis can be calibrated in terms of the strength of standard theories \cite{keisler:06}.
\item The axiom of choice and the continuum hypothesis are conservative extensions
  of set theory for $\Sigma^2_1$ sentences in the analytic hierarchy \cite{kunen:80}.
\end{enumerate}
Such results draw on a variety of methods. Some can be obtained by direct translation of one theory into another. Many are proved using cut-elimination or normalization \cite{buss:98b,sieg:85}. The double-negation translation is a remarkably effective tool when it comes to reducing classical theories to constructive ones, and can often be supplemented by realizability, functional interpretation, or other arguments \cite{avigad:00b,kohlenbach:08,troelstra:98}. Model-theoretic methods can often be used, though they do not provide specific algorithms to carry out the translation \cite{hajek:pudlak:93,kaye:91}. Even forcing methods, originally developed as a set-theoretic technique, can be fruitfully be applied in proof-theoretic settings \cite{avigad:04,kunen:80}. 

\subsection{Characterizing logical strength}

The results described in the previous section serve to characterize the strength of one axiomatic theory in terms of another. Showing that a theory $T_2$ is conservative over $T_1$ shows that, in particular, $T_2$ is consistent, if $T_1$ is. This provides a comparison of the \emph{consistency strength} of the two theories. 

But there are other ways of characterizing the strength of a theory. For example, the notion of an \emph{ordinal} generalizes the notion of a counting number. Starting with the natural numbers, we can add an infinite ``number,'' $\omega$, and keep going:
\[
0, 1, 2, 3, \ldots, \omega, \omega + 1, \omega + 2, \omega + 3, \ldots
\]
We can then proceed to add even more exotic numbers, like $\omega \cdot 2$, $\omega^2$, and $\omega^\omega$. 
The ordering on these particular expressions is computable, in the sense that one can write a computer program to compare any two them. What makes them ordinals is that they satisfy a principle of \emph{transfinite induction,} which generalizes the principle of induction on the natural numbers. Ordinal analysis gauges the strength of a theory in terms of such computable ordinals: the stronger a theory is, the more powerful the principles of transfinite induction it can prove. See, for example, \cite{pohlers:98,pohlers:09,takeuti:87}.

Alternatively, one can focus on a theory's \emph{computational strength}. Suppose a theory $T$ proves a statement of the form $\fa x \ex y R(x,y)$, where $x$ and $y$ range over the natural numbers, and $R$ is a computationally decidable predicate. This tells us that a computer program that, on input $x$, systematically searches for a $y$ satisfying $R(x,y)$ always succeeds in finding one. Now suppose $f$ is a function that, on input $x$, returns a value that is easily computed from the least $y$ satisfying $R(x,y)$. For example, $R(x,y)$ may assert that $y$ codes a halting computation of a particular Turing machine on input $x$, and $f$ may return the result of such a computation. Then $f$ is a computable function, and we can say that the theory, $T$, proves that $f$ is totally defined on the natural numbers. A simple diagonalization shows that no effectively axiomatized theory can prove the totality of every computable function in this way, so this suggests using the set of computable functions that the theory can prove to be total as a measure of its strength.

A number of theories have been analyzed in these terms. For example, by the results in the last section, the provably total computable functions of $\na{PRA}$, $\na{I\Sigma_1}$, $\na{RCA_0}$, and $\na{WKL_0}$ are all the primitive recursive functions. In contrast, one can characterize the provably total computable functions of $\na{PA}$ and $\na{HA}$ in terms of higher-type primitive recursion \cite{avigad:feferman:98,troelstra:98}, or using principles of primitive recursion along an ordinal known as $\varepsilon_0$ \cite{takeuti:87,pohlers:09}. Weaker theories of arithmetic can be used to characterize complexity classes like the polynomial time computable functions \cite{buss:98b}.

\section{Applications}
\label{applications:section}

In this final section, I will describe some of the ways that proof theory interacts with other disciplines. As emphasized in Section~\ref{introduction:section}, I am only considering applications of the traditional, metamathematical branch of proof theory. Formal deductive methods, more broadly, have applications across philosophy and the sciences, and the use of proof-theoretic methods in the study of these formal deductive systems is far too diverse to survey here.

\subsection{Proof mining}

One way in which traditional proof-theoretic methods have been applied is in the process of extracting useful information from ordinary mathematical proofs. The reductive results of the twentieth century showed, in principle, that many classical proofs can be interpreted in constructive terms. In practice, these ideas have been adapted and extended to the analysis of ordinary mathematical proofs. Georg Kreisel described the process of extracting such information as ``unwinding proofs,'' and Ulrich Kohlenbach has more recently adopted the name ``proof mining'' \cite{kohlenbach:08}.

Substantial work is needed to turn this vague idea into something practicable. Ordinary mathematical proofs are not presented in formal systems, so there are choices to be made in the formal modeling. In addition, the general metamathematical tools have to be tailored and adjusted to yield the information that is sought in particular domains. Thus the work requires a deep understanding of both the proof-theoretic methods and the domain of mathematics in question. The field has already had a number of successes in fields like functional analysis and ergodic theory; see, for example, \cite{kohlenbach:08}.

\subsection{Combinatorial independences}

Yet another domain where a syntactic, foundational perspective is important is in the search for natural combinatorial independences, that is, natural finitary combinatorial principles that are independent of
conventional mathematical methods. The Paris-Harrington statement \cite{paris:harrington:77} is an early example of such a principle. Since then, Harvey Friedman, in particular, has long sought to find exotic combinatorial behavior in familiar mathematical settings. Such work gives us glimpes into what goes on just beyond ordinary patterns of mathematical reasoning, and yields interesting mathematics as well. See the extensive introduction to \cite{friedman:unp:b} for an overview of results in this area.

\subsection{Constructive mathematics and type theory}

As noted above, proof theory is often linked with constructive mathematics, for historical reasons. After all, Hilbert's program was initially an attempt to justify mathematics with respect to methods that are finitary, which is to say, syntactic, algorithmic, and impeccably constructive. Contemporary work in constructive mathematics and type theory draws on the following facts:
\begin{itemize}
\item Logical constructions can often be interpreted as programming principles.
\item Conversely, programming principles can be interpreted as logical constructions.
\item One can thereby design (constructive) proof systems that combine aspects of both programming and proving.
\end{itemize}
The references in Section~\ref{type:theory} above provide logical perspectives on constructive type theory. For a computational  perspective, see \cite{pierce:04}.

\subsection{Automated reasoning and formal verification}

Another domain where proof-theoretic methods are of central importance is in the field of automated reasoning and formal verification. In computer science, researchers use formal methods to help verify that hardware and software are bug-free and conform to their specifications. Moreover, recent developments have shown that computational formal methods can be used to help verify the correctness of complex mathematical proofs as well. Both efforts have led to interactive approaches, whereby a user works with a computational proof assistant to construct a formal proof of the relevant claims. The have also led to more automated approaches, where software is supposed to carry out the task with little user input. In both cases, proof-theoretic methods are invaluable, for designing the relevant logical calculi, for isolating features of proofs that enable one to cut down the search space and traverse it effectively, and for replacing proof search with calculation wherever possible. 

For more information on automated reasoning, see \cite{harrison:09b,robinson:voronkov:01}. For more information on formally verified mathematics, see \cite{wiedijk:06}, or the the December 2008 issue of the \emph{Notices of the American Mathematical Society}, which was devoted to formal proof.

\subsection{Proof complexity}

Finally, the field of proof complexity combines methods and insights from proof theory and computational complexity. For example, the complexity class NP can be viewed as the class of problems for which an affirmative answer has a short (polynomial-size) proof in a suitable calculus. Thus the conjecture that NP is not equal to co-NP (which is weaker than saying P is not equal to NP) is equivalent to saying that in general there is no propositional calculus that has efficient proofs of every tautology. Stephen Cook has suggested that one way of building up to the problem is to show that \emph{particular} proof systems are not efficient, by establishing explicit lower bounds. Such information is also of interest in automated reasoning, where one wishes to have a detailed understanding of the types of problems that can be expected to have short proofs in various calculi. The works \cite{krajicek:95,pudlak:98,segerlind:07,urquhart:95} provide excellent introductory overviews.


\end{document}